# EXTREME VALUE THEORY, ERGODIC THEORY AND THE BOUNDARY BETWEEN SHORT MEMORY AND LONG MEMORY FOR STATIONARY STABLE PROCESSES[1]

BY GENNADY SAMORODNITSKY

*Cornell University*

We study the partial maxima of stationary $\alpha$-stable processes. We relate their asymptotic behavior to the ergodic theoretical properties of the flow. We observe a sharp change in the asymptotic behavior of the sequence of partial maxima as flow changes from being dissipative to being conservative, and argue that this may indicate a change from a short memory process to a long memory process.

**1. Introduction.** Let $\mathbf{X} = (X_0, X_1, \ldots)$ be a stationary symmetric $\alpha$-stable (S$\alpha$S) process, $0 < \alpha < 2$. How does one decide whether or not the process has long range dependence?

Since $\alpha$-stable random variables with $0 < \alpha < 2$ have infinite second moment, one cannot use correlations to tell when a stationary $\alpha$-stable process has long range dependence. Covariance-like functions have been tried [see, e.g., Astrauskas, Levy and Taqqu (1991)], but their usefulness seems to be limited. In fact, even for stationary processes with a finite second moment, the definition of long range dependence based on the true correlation function is of uncertain value unless the process is a Gaussian process, or very close to being one.

Instead of using the correlation function or looking for a substitute, we propose a different approach. Suppose that $(\mathcal{P}_\theta, \theta \in \Theta)$ is a family of laws of a stationary stochastic process $(X_0, X_1, X_2, \ldots)$, where $\Theta$ is some parameter space. Assume that the marginal laws of the process do not change much as $\theta$ varies (perhaps, the marginal laws remain constant, or only the global scale changes, if we are considering, say, Gaussian or S$\alpha$S processes). Suppose we

Received July 2002; revised April 2003.
[1]Supported in part by NSF Grant DMS-00-71073 at Cornell University.

are given a functional of interest $R$, a (measurable) functional on $\mathbb{R}^\infty$. The behavior of this functional is different, in general, under different laws $\mathcal{P}_\theta$. Suppose that there is a partition of the parameter space $\Theta$ into two parts, $\Theta_0$ and $\Theta_1$, such that the behavior of the functional changes dramatically as one crosses the boundary between $\Theta_0$ and $\Theta_1$. Such change may be caused by various factors (e.g., by changing heaviness of the tails), but in some cases it may make sense to talk about that boundary as a boundary between short range dependence and long range dependence. That is, the change from short memory to long memory occurs as a phase transition. We emphasize that the behavior of each individual functional does not define short or long memory, and the phase transition indicated should occur for a large group of functionals for a boundary to be called a change from short to long memory. A complete theory is missing at the moment. In this paper we find one important functional undergoing such a phase transition at a boundary.

Existence of boundaries with such properties has been observed before. For stationary zero mean Gaussian processes parameterized by the common variance, and by the correlation function, such a phase transition occurs when the correlations stop being summable. The functional of interest here is the sequence of the partial sums, and its distributional rate of growth changes significantly at the boundary. The rate of growth of the partial sums may change its order of magnitude whether or not the second moment is finite. This has been observed many times on the example of the increments of self-similar processes with stationary increments. In the Gaussian case the family of such processes are fractional Brownian motions, parameter $H$ of self-similarity has to be in the interval $(0, 1)$, and the partial sums of the increment process (the so-called fractional Gaussian noise) increase at the rate higher than $n^{1/2}$ when $H > 1/2$. Hence $H = 1/2$ is considered to be the boundary between short and long memory for fractional Gaussian noise. See, for example, Mandelbrot (1975) and Mandelbrot and Taqqu (1979), as well as a more recent discussion in Beran (1994). A similar phenomenon occurs for the increments of self-similar $\alpha$-stable processes with stationary increments, $1 < \alpha < 2$, that have infinite variance. Here the range of parameter $H$ of self-similarity is still $(0, 1)$, and the boundary where the partial sums of the increment process start increasing at the rate higher than the i.i.d. case (i.e., faster than $n^{1/\alpha}$) is that of $H = 1/\alpha$. See, for example, Samorodnitsky and Taqqu (1994). No such boundary is possible if we consider the increments of self-similar $\alpha$-stable processes with stationary increments and $0 < \alpha \leq 1$. If one uses the boundary $H = 1/\alpha$ to define long range dependence, one would have to conclude that long range dependence is impossible if $0 < \alpha \leq 1$.

The present paper uncovers a different boundary for stationary $\alpha$-stable processes. We believe that this is a very fundamental boundary, and it is based on ergodic-theoretical properties of nonsingular flows underlying such processes. Specifically, we concentrate on the partial maxima sequence

(1.1) $$M_n = \max(|X_0|, |X_1|, \ldots, |X_{n-1}|), \qquad n = 1, 2, \ldots,$$

and its distributional rate of growth. We will see that the parameter space consists of two parts, in one of which the partial maxima grow at the rate $n^{1/\alpha}$, which is the rate at which partial maxima of i.i.d. $\alpha$-stable random variables grow, while in the other part of the parameter space the partial maxima grow at a strictly slower rate. Moreover, in the latter part of the parameter space the actual rate of growth may depend on the choice of the parameters. This boundary is present for all $0 < \alpha < 2$. Again, by itself the change in behavior of a single functional does not qualify this boundary as that between short and long memory. Our conjecture is, however, that many other important changes occur at that boundary.

This paper is organized as follows. In Section 2 we provide a background on integral representations of stationary S$\alpha$S processes, and elements of the theory developed by Rosiński (1994, 1995) relating such integral representations and ergodic theory. In Section 3 we discuss the behavior of a certain deterministic sequence controlling the rate of growth of the partial maxima. In Section 4 we prove the main result, Theorem 4.1, dealing with the asymptotic behavior of the sequence partial maxima. Section 5 discusses a number of examples illustrating the results of Section 4. A brief conclusion discussing what the results of this paper tell us about short and long memory for stationary stable processes is in Section 6.

**2. Ergodic theory and representations of stationary $\alpha$-stable processes.** Throughout this section $\mathbf{X} = (X_0, X_1, \ldots)$ is a S$\alpha$S process, $0 < \alpha < 2$. Every (not necessarily stationary) S$\alpha$S process has an integral representation

(2.1) $$X_n = \int_E f_n(x) M(dx), \qquad n = 0, 1, 2, \ldots,$$

where $M$ is a S$\alpha$S random measure on a measurable space $(E, \mathcal{E})$ with a $\sigma$-finite control measure $m$, while $f_n \in L^\alpha(m, \mathcal{E})$ for all $n$. See Chapter 3 in Samorodnitsky and Taqqu (1994) on $\alpha$-stable random measures and integrals with respect to these measures, and Chapter 13 there on integral representations as above as well as on the history of such representations.

When a process $\mathbf{X}$ is stationary, the integral representation can be selected to be of a particular form, according to a theory developed in a series of papers by Rosiński. See, for example, Rosiński (1995); various facts presented below can be found in that paper. Specifically, a stationary S$\alpha$S process has an integral representation of the form (2.1) with

(2.2) $$f_n(x) = a_n(x) \left( \frac{dm \circ \phi^n}{dm}(x) \right)^{1/\alpha} f \circ \phi^n(x), \qquad x \in E,$$

for $n = 0, 1, 2, \ldots$, where $\phi : E \to E$ is a measurable nonsingular map (i.e., a one-to-one map with both $\phi$ and $\phi^{-1}$ measurable, mapping the control measure $m$ into an equivalent measure),

$$a_n(x) = \prod_{j=0}^{n-1} u \circ \phi^j(x), \qquad x \in E,$$

for $n = 0, 1, 2, \ldots$, with $u : E \to \{-1, 1\}$ a measurable function and $f \in L^\alpha(m, \mathcal{E})$.

That is, the process **X** is determined by a single function $f \in L^\alpha(m, \mathcal{E})$, a *cocycle* $(a_n, n = 0, 1, 2, \ldots)$ and a *flow* $(\phi^n, n = 0, 1, 2, \ldots)$. This triple [taken together with the space $(E, \mathcal{E}, m)$ on which it lives] can, therefore, be taken as a parameterization of stationary S$\alpha$S processes. When working with this parameterization, the task is to relate the ergodic-theoretic properties of the flow to the probabilistic properties of the stable process. We are interested in properties that we can interpret as related to the length of memory of a stationary stable process.

Let $E = C \cup D$ be the Hopf decomposition of $E$ with respect to the flow. That is, $C$ and $D$ are measurable $\phi$-invariant sets, such that the flow is conservative on $C$ and dissipative on $D$; see Krengel (1985) for the various ergodic-theoretical notions and facts we use in this paper. Writing

$$
\begin{aligned}
X_n &= \int_C a_n(x) \left( \frac{dm \circ \phi^n}{dm}(x) \right)^{1/\alpha} f \circ \phi^n(x) M(dx) \\
&\quad + \int_D a_n(x) \left( \frac{dm \circ \phi^n}{dm}(x) \right)^{1/\alpha} f \circ \phi^n(x) M(dx) \\
&=: X_n^C + X_n^D, \qquad n = 0, 1, 2, \ldots,
\end{aligned}
\tag{2.3}
$$

leads to a unique in law decomposition of a stationary S$\alpha$S process into a sum of two independent such processes, one of which is generated by a conservative flow, and the other by a dissipative flow.

Intuitively, one expects stable processes generated by conservative flows to have a longer memory than those generated by dissipative flows, simply because a conservative flow "tends to keep coming back," and so the same values of the random measure $M$ contribute to observations $X_n$ far separated in time. Consider, for example, stationary processes generated by a dissipative flow. Such a process has a *mixed moving average* representation of the form

$$X_n = \int_W \int_{\mathbb{Z}} f(v, x - n) M(dv, dx), \qquad n = 0, 1, 2, \ldots, \tag{2.4}$$

with $M$ a S$\alpha$S random measure on a product measurable space $(W \times \mathbb{Z}, \mathcal{W} \times \mathcal{B})$ with control measure $m = \nu \times l$, where $\nu$ is a $\sigma$-finite measure on $(W, \mathcal{W})$, $l$ is the counting measure on $\mathbb{Z}$ and $f \in L^\alpha(m, \mathcal{W} \times \mathcal{B})$. Such processes are always mixing irrespective of what the kernel $f$ in either (2.2) or (2.4) is [or what the cocycle in (2.2) is]. See Surgailis, Rosiński, Mandrekar and Cambanis (1993). On the other hand, stable processes generated by conservative flows are often not even ergodic. For example, if a conservative flow is measure preserving and the expected return time to any set of a finite positive measure is finite, then the stable process is not ergodic [see Rosiński and Samorodnitsky (1996)]. Another example demonstrating that conservative flows tend to lead to a longer memory (in the case $1 < \alpha < 2$) can be found in Mikosch and Samorodnitsky (2000), who studied ruin probabilities.

**3. The sequence $(b_n)$.** Let $\mathbf{X} = (X_0, X_1, \dots)$ be a stationary S$\alpha$S process, $0 < \alpha < 2$. We assume that the process is given in an integral representation of the form (2.1), where $(f_n)$ is of the form (2.2). It turns out that, to a large extent, the asymptotic behavior of the maximal functional $M_n$ in (1.1) is related to the quantity defined below. Let

$$(3.1) \qquad b_n = \left( \int_E \max_{j=0,\dots,n-1} |f_j(x)|^\alpha m(dx) \right)^{1/\alpha}, \qquad n = 1, 2, \dots.$$

In fact, to a certain extent $b_n$ controls "the size" of $M_n$ even without the assumption of stationarity of the process. Indeed, for any $0 < p < \alpha$, there are constants $c_{\alpha,p}, C_{\alpha,p} \in (0, \infty)$ such that, for $1 < \alpha < 2$,

$$(3.2) \qquad c_{\alpha,p} \leq \frac{1}{b_n}(EM_n^p)^{1/p} \leq C_{\alpha,p}(\log n)^{1/\alpha'},$$

where $\alpha'$ is the conjugate of $\alpha$ in $1/\alpha + 1/\alpha' = 1$, while for $\alpha = 1$,

$$(3.3) \qquad c_{1,p} \leq \frac{1}{b_n}(EM_n^p)^{1/p} \leq C_{1,p} L_2 n,$$

where $L_2 n = \max(1, \log \log n)$. Finally, for $0 < \alpha < 1$,

$$(3.4) \qquad c_{\alpha,p} \leq \frac{1}{b_n}(EM_n^p)^{1/p} \leq C_{\alpha,p}.$$

See Theorem 2.1 in Marcus (1984).

We will see that for stationary S$\alpha$S processes the sequence $(b_n)$ tells us even more about the sequence $(M_n)$ of the partial maxima. Note that the sequence $(b_n)$ is completely determined by the process, and does not depend on a particular integral representation. Certain important features of this sequence are determined by the flow in (2.2) underlying the process. In particular, the next result shows that the sequence $b_n$ grows at slower rate for processes generated by a conservative flow than for processes generated by a dissipative flow.

THEOREM 3.1. *Let $(f_n)$ be given by* (2.2).

(i) *If the flow $(\phi^n)$ is conservative, then:*

$$n^{-1/\alpha} b_n \to 0 \qquad \text{as } n \to \infty. \tag{3.5}$$

(ii) *If the flow is dissipative, and the process is given in the mixed moving average form* (2.4), *then:*

$$\lim_{n \to \infty} n^{-1/\alpha} b_n = \left( \int_W g(v)^\alpha \nu(dv) \right)^{1/\alpha} \in (0, \infty), \tag{3.6}$$

*where*

$$g(v) = \sup_{k=0, \pm 1, \pm 2, \ldots} |f(v, k)| \qquad \text{for } v \in W. \tag{3.7}$$

PROOF. (i) Suppose first that the flow $(\phi^n)$ is control measure $m$ preserving, and that $f = \mathbf{1}_A$ for some $A \in \mathcal{E}$ with $0 < m(A) < \infty$. In that case the Radon–Nykodim derivative in (2.2) disappears, and we get

$$b_n^\alpha = m\left( \bigcup_{k=0}^{n-1} \phi^{-k}(A) \right), \qquad n = 1, 2, \ldots. \tag{3.8}$$

We use a construction used to prove the Kac recurrence theorem; see, for example, Theorem 3.6 in Krengel (1985). Let

$$\tau_A(x) = \inf\{k \geq 1 : \phi^k(x) \in A\} \qquad \text{for } x \in A$$

be the first recurrence time to $A$, and let, for $k = 1, 2, \ldots$,

$$R_k = \{x \in A : \tau_A(x) = k\}$$

and

$$A_k = \{x \in A^c : \phi^k(x) \in A, \phi^j(x) \notin A, j = 1, \ldots, k-1\}$$

be, correspondingly, the set of points in $A$ returning to $A$ after $k$ steps, and the set of points outside of $A$ entering $A$ for the first time after $k$ steps. Let also $A_0 = A$. Note that

$$m\left( \bigcup_{k=0}^{n-1} \phi^{-k}(A) \right) = \sum_{k=0}^{n-1} m(A_k). \tag{3.9}$$

Furthermore, for every $k = 1, 2, \ldots$,

$$\begin{aligned} R_k &= \{x \in E : \phi^k(x) \in A, \phi^j(x) \notin A, j = 1, \ldots, k-1\} \\ &\quad - \{x \in A^c : \phi^k(x) \in A, \phi^j(x) \notin A, j = 1, \ldots, k-1\} \\ &= \phi^{-1}(A_{k-1}) - A_k, \end{aligned}$$

and so by the measure-preserving property of the flow,

$$m(R_k) = m(A_{k-1}) - m(A_k).$$

EXTREMES OF STATIONARY STABLE PROCESSES 7Summing up, we see that, for every $k = 0, 1, \ldots,$

$$(3.10) \qquad m(A_k) = \sum_{j=k+1}^{\infty} m(R_j) + \lim_{n \to \infty} m(A_n).$$

Since the flow is conservative, we can use (3.10) with $k = 0$ to see that the limit above is equal to zero. Therefore, by (3.9) we conclude that

$$\frac{1}{n} m\left(\bigcup_{k=0}^{n-1} \phi^{-k}(A)\right) \to 0 \qquad \text{as } n \to \infty,$$

which establishes (3.5) for indicator functions in the measure-preserving case.

To establish (3.5) in the general case (i.e., not necessarily measure-preserving conservative flow and a general kernel $f$), denote

$$w_n(x) = \frac{dm \circ \phi^n}{dm}(x), \qquad n = 0, 1, 2, \ldots, \ x \in E.$$

Then

$$\phi_*(x, y) = \left(\phi(x), \frac{y}{w_1(x)}\right), \qquad x \in E, \ y > 0,$$

defines a conservative flow on $(E \times (0, \infty), \mathcal{E} \times \mathcal{B}, m \times \text{Leb})$ which can be expressed as

$$\phi_*^n(x, y) = \left(\phi^n(x), \frac{y}{w_n(x)}\right), \qquad x \in E, \ y > 0 \text{ for } n = 0, \pm 1, \pm 2, \ldots,$$

and this flow preserves the measure $m \times \text{Leb}$; see Maharam (1964).

Let

$$A = \{(x, y) \in E \times (0, \infty) : 0 < y \leq |f(x)|^\alpha\}.$$

Note that

$$m \times \text{Leb}(A) = \int_E |f(x)|^\alpha \, m(dx) \in (0, \infty).$$

Furthermore,

$$m \times \text{Leb}\left(\bigcup_{k=0}^{n-1} \phi_*^{-k}(A)\right)$$

$$= \int_E \int_0^\infty \mathbf{1}_{\bigcup_{k=0}^{n-1} \phi_*^{-k}(A)}(x, y) \, m(dx) \, dy$$

$$= \int_E \int_0^\infty \max_{k=0,\ldots,n-1} \mathbf{1}_{\phi_*^{-k}(A)}(x, y) \, m(dx) \, dy$$

$$(3.11) \quad = \int_E \int_0^\infty \max_{k=0,\ldots,n-1} \mathbf{1}(0 < y \le w_k(x)|f \circ \phi^k(x)|^\alpha) \, m(dx) \, dy$$

$$= \int_E \max_{k=0,\ldots,n-1} w_k(x)|f \circ \phi^k(x)|^\alpha \, m(dx)$$

$$= \int_E \max_{k=0,\ldots,n-1} |f_k(x)|^\alpha \, m(dx) = b_n^\alpha.$$

Since the left-hand side of (3.11) is $o(n)$ by the already considered case of measure-preserving flows and indicator functions, this establishes (3.5) in full generality.

(ii) We start with the case where $f$ has a compact support, that is,

$$(3.12) \quad f(v,k)\mathbf{1}_{W\times[-m,m]^c}(v,k) \equiv 0 \quad \text{for some } m=1,2,\ldots.$$

In that case, for $n \ge 2m+1$, we have

$$b_n^\alpha = \sum_{j=-m-n+1}^{m} \int_W \max_{k=0,\ldots,n-1} |f(v,j+k)|^\alpha \, \nu(dv)$$

$$= \sum_{j=m-n+1}^{-m} \int_W \max_{k=0,\ldots,n-1} |f(v,j+k)|^\alpha \, \nu(dv)$$

$$+ \sum_{j=-m-n+1}^{m-n} \int_W \max_{k=0,\ldots,n-1} |f(v,j+k)|^\alpha \, \nu(dv)$$

$$+ \sum_{j=-m+1}^{m} \int_W \max_{k=0,\ldots,n-1} |f(v,j+k)|^\alpha \, \nu(dv) =: T_n + R_n^{(1)} + R_n^{(2)}.$$

Observe that, for each $j = m - n + 1, \ldots, -m$,

$$\max_{k=0,\ldots,n-1} |f(v,j+k)| = g(v),$$

while

$$\max_{k=0,\ldots,n-1} |f(v,j+k)| \le g(v)$$

for other values of $j$. Therefore,

$$T_n = (n-2m) \int_W g(v)^\alpha \, \nu(dv),$$

while

$$|R_n^{(i)}| \le 2m \int_W g(v)^\alpha \, \nu(dv) \quad \text{for } i = 1, 2.$$

Therefore (3.6) in the case of a compact supported $f$ follows. In the general case, given $\epsilon > 0$, choose a compact supported $f_\epsilon$ such that $|f_\epsilon(v,k)| \leq |f(v,k)|$ for all $v, k$, and
$$\sum_{k=-\infty}^{\infty} \int_W |f(v,k)|^\alpha \nu(dv) - \sum_{k=-\infty}^{\infty} \int_W |f_\epsilon(v,k)|^\alpha \nu(dv) \leq \epsilon.$$
Let
$$g_\epsilon(v) = \sup_{k=0,\pm 1, \pm 2, \ldots} |f_\epsilon(v,k)|,$$
$v \in W$. Then
$$0 \leq \int_W g(v)^\alpha \nu(dv) - \int_W g_\epsilon(v)^\alpha \nu(dv)$$
$$\leq \int_W \sup_{k=0,\pm 1, \pm 2, \ldots} (|f(v,k)|^\alpha - |f_\epsilon(v,k)|^\alpha) \nu(dv)$$
$$\leq \int_W \sum_{k-\infty}^{\infty} (|f(v,k)|^\alpha - |f_\epsilon(v,k)|^\alpha) \nu(dv)$$
$$= \sum_{k=-\infty}^{\infty} \int_W |f(v,k)|^\alpha \nu(dv) - \sum_{k=-\infty}^{\infty} \int_W |f_\epsilon(v,k)|^\alpha \nu(dv) \leq \epsilon.$$
Therefore,
$$\left| \frac{1}{n} b_n^\alpha - \int_W g(v)^\alpha \nu(dv) \right|$$
$$\leq \frac{1}{n} \left| \sum_{j=-\infty}^{\infty} \int_W \max_{k=0,\ldots,n-1} |f(v, j+k)|^\alpha \nu(dv) \right.$$
$$\left. - \sum_{j=-\infty}^{\infty} \int_W \max_{k=0,\ldots,n-1} |f_\epsilon(v, j+k)|^\alpha \nu(dv) \right|$$
$$+ \left| \frac{1}{n} \sum_{j=-\infty}^{\infty} \int_W \max_{k=0,\ldots,n-1} |f_\epsilon(v,j+k)|^\alpha \nu(dv) - \int_W g_\epsilon(v)^\alpha \nu(dv) \right|$$
$$+ \left| \int_W g_\epsilon(v)^\alpha \nu(dv) - \int_W g(v)^\alpha \nu(dv) \right| =: T_n^{(1)} + T_n^{(2)} + T_n^{(3)}.$$
By the above, $T_n^{(3)} \leq \epsilon$, and the same argument shows that $T_n^{(1)} \leq \epsilon$ as well. Furthermore, by the already considered compact support case, $T_n^{(2)} \to 0$ as $n \to \infty$. Hence
$$\limsup_{n \to \infty} \left| \frac{1}{n} b_n^\alpha - \int_W g(v)^\alpha \nu(dv) \right| \leq 2\epsilon,$$
and, since $\epsilon > 0$ is arbitrary, the proof of (3.6) is complete. $\square$

**4. Maxima of stationary stable processes.** In this section we investigate the rate of growth of the sequence $(M_n)$ of partial maxima of a stationary S$\alpha$S process, $0 < \alpha < 2$. We will see, in particular, that if such a process has a nonzero component $\mathbf{X}^C$ in (2.3) generated by a dissipative flow, then the partial maxima grow at the rate of $n^{1/\alpha}$, while if the process is generated by a conservative flow, then the partial maxima grow at a slower rate. The following is the main theorem of this paper.

THEOREM 4.1. *Let $\mathbf{X} = (X_0, X_1, \ldots)$ be a stationary S$\alpha$S process, $0 < \alpha < 2$, with integral representation* (2.1), *and functions* $(f_n)$ *given by* (2.2).

(i) *Suppose that $\mathbf{X}$ is not generated by a conservative flow [i.e., that the component $\mathbf{X}^D$ in* (2.3) *generated by a dissipative flow is nonzero]. Then*

$$(4.1) \qquad \frac{1}{n^{1/\alpha}} M_n \Rightarrow C_\alpha^{1/\alpha} K_X Z_\alpha$$

*weakly as $n \to \infty$, where*

$$K_X = \left( \int_W g(v)^\alpha \nu(dv) \right)^{1/\alpha}$$

*and $g$ given by* (3.7) *for any representation of $\mathbf{X}^D$ in the mixed moving average form* (2.4). *Furthermore,*

$$(4.2) \quad C_\alpha = \left( \int_0^\infty x^{-\alpha} \sin x\, dx \right)^{-1} = \begin{cases} \dfrac{1-\alpha}{\Gamma(2-\alpha)\cos(\pi\alpha/2)}, & \text{if } \alpha \neq 1, \\ \dfrac{2}{\pi}, & \text{if } \alpha = 1, \end{cases}$$

*and $Z_\alpha$ is the standard Frechét-type extreme value random variable with the distribution*

$$P(Z_\alpha \leq z) = e^{-z^{-\alpha}}, \qquad z > 0.$$

(ii) *Suppose that $\mathbf{X}$ is generated by a conservative flow. Then*

$$(4.3) \qquad \frac{1}{n^{1/\alpha}} M_n \to 0$$

*in probability as $n \to \infty$. Furthermore, with $b_n$ given by* (3.1),

$$(4.4) \quad \left( \frac{1}{c_n} M_n \right) \qquad \text{is not tight for any positive sequence } c_n = o(b_n),$$

*while*

$$(4.5) \left( \frac{1}{b_n f_n} M_n \right) \quad \text{is tight, where} \quad f_n = \begin{cases} (\log n)^{1/\alpha'}, & \text{if } 1 < \alpha < 2, \\ L_2 n, & \text{if } \alpha = 1, \\ 1, & \text{if } 0 < \alpha < 1. \end{cases}$$

EXTREMES OF STATIONARY STABLE PROCESSES        11*If, for some $\theta > 0$ and $c > 0$,*

(4.6) $$b_n \geq cn^\theta \qquad \text{for all } n \geq 1,$$

*then (4.5) holds with $f_n \equiv 1$ for all $0 < \alpha < 2$.*

*Furthermore, for $n = 1, 2, \ldots$, let $\eta_n$ be a probability measure on $(E, \mathcal{E})$ with*

(4.7) $$\frac{d\eta_n}{dm}(x) = b_n^{-\alpha} \max_{j=0,\ldots,n-1} |f_j(x)|^\alpha, \qquad x \in E,$$

*and let $U_j^{(n)}$, $j = 1, 2$, be independent $E$-valued random variables with common law $\eta_n$. Suppose that (4.6) holds and that, in addition, for any $\epsilon > 0$,*

(4.8) $$P\bigg(\text{for some } k = 0, 1, \ldots, n-1,$$
$$\frac{|f_k(U_j^{(n)})|}{\max_{i=0,\ldots,n-1}|f_i(U_j^{(n)})|} > \epsilon, \ j = 1, 2\bigg) \to 0$$

*as $n \to \infty$. Then*

(4.9) $$\frac{1}{b_n} M_n \Rightarrow C_\alpha^{1/\alpha} Z_\alpha$$

*weakly as $n \to \infty$.*

REMARK 4.2. Here are some sufficient conditions for (4.8). If

(4.10) $$\lim_{n \to \infty} \frac{b_n}{n^{1/2\alpha}} = \infty,$$

then (4.8) holds. Indeed, let $r_n$ denote the probability in the left-hand side of (4.8). Clearly,

$$r_n \leq \sum_{k=0}^{n-1} \bigg(P\bigg(\frac{|f_k(U_1^{(n)})|}{\max_{i=0,\ldots,n-1}|f_i(U_1^{(n)})|} > \epsilon\bigg)\bigg)^2.$$

Furthermore, for every $k = 0, 1, \ldots, n-1$,

$$P\bigg(\frac{|f_k(U_1^{(n)})|}{\max_{i=0,\ldots,n-1}|f_i(U_1^{(n)})|} > \epsilon\bigg)$$
$$= b_n^{-\alpha} \int_E \mathbf{1}\bigg(\frac{|f_k(x)|}{\max_{i=0,\ldots,n-1}|f_i(x)|} > \epsilon\bigg) \max_{i=0,\ldots,n-1} |f_i(x)|^\alpha m(dx)$$
$$\leq \epsilon^{-\alpha} b_n^{-\alpha} \int_E |f_k(x)|^\alpha m(dx),$$

and (4.8) follows from (4.10) since, by the stationarity, the last integral does not depend on $k$.

Alternatively, assume that $m$ is a finite measure, the flow is measure $m$ preserving, the sequence $(b_n^{-\alpha} \max_{j=0,\ldots,n-1} |f_j(x)|^\alpha)$, $n = 1, 2, \ldots$, is uniformly integrable with respect to $m$ and, for every $\epsilon > 0$

$$(4.11) \qquad \lim_{n \to \infty} n^{1/2} m\{x \in E : |f(x)| > \epsilon b_n\} = 0,$$

where $f$ is the kernel in (2.2). Then (4.8) holds.

Indeed, let $\|m\|$ be the total mass of $m$. Given a $\delta > 0$, select $M > 0$ such that

$$\int_E \mathbf{1}\left(\max_{i=0,\ldots,n-1} |f_i(x)|^\alpha > M b_n^\alpha\right) \max_{i=0,\ldots,n-1} |f_i(x)|^\alpha \, m(dx) \leq \delta \, b_n^\alpha$$

for all $n \geq 1$. We have with $\epsilon$ from (4.8),

$$r_n \leq 4\delta + b_n^{-2\alpha} \sum_{k=0}^{n-1} \bigg(\int_E \max_{j=0,\ldots,n-1} |f_j(x)|^\alpha$$
$$\times \mathbf{1}\bigg(\delta \|m\|^{-1} b_n^\alpha \leq \max_{j=0,\ldots,n-1} |f_j(x)|^\alpha \leq M b_n^\alpha\bigg)$$
$$\times \mathbf{1}\bigg(\frac{|f_k(x)|}{\max_{i=0,\ldots,n-1} |f_i(x)|} > \epsilon\bigg) m(dx)\bigg)^2$$
$$\leq 4\delta + M^2 n \bigg(\int_E \mathbf{1}(|f(x)| > \epsilon \delta \|m\|^{-1} b_n) m(dx)\bigg)^2.$$

Therefore, using (4.11), we obtain

$$\limsup_{n \to \infty} r_n \leq 4\delta,$$

and (4.8) follows by letting $\delta \to 0$.

PROOF OF THEOREM 4.1. We use a series representation of the random vector $(X_0, X_1, \ldots, X_{n-1})$ of the form

$$(4.12) \qquad X_k = b_n C_\alpha^{1/\alpha} \sum_{j=1}^\infty \varepsilon_j \Gamma_j^{-1/\alpha} \frac{f_k(U_j^{(n)})}{\max_{i=0,\ldots,n-1} |f_i(U_j^{(n)})|},$$
$$k = 0, 1, \ldots, n-1,$$

where $C_\alpha$ is given by (4.2), $\varepsilon_1, \varepsilon_2, \ldots$ are i.i.d. Rademacher random variables (symmetric $\pm 1$-valued random variables), $\Gamma_1, \Gamma_2, \ldots$ is a sequence of the arrival times of a unit rate Poisson process on $(0, \infty)$, and $(U_j^{(n)})$ are i.i.d. $E$-valued random variables with common law given by (4.7). All three sequences are independent. See Section 3.10 in Samorodnitsky and Taqqu (1994). Of course, the representation in (4.12) is in law.

We start with observing that (4.5) follows from (3.2)–(3.4) regardless of the properties of the flow. To check (4.4), we use the series representation (4.12) and symmetry. Let $\mathcal{G}$ be the $\sigma$-field generated by $\varepsilon_1, (\Gamma_j, j \geq 1)$ and $(U_j^{(n)}, j \geq 1)$. Letting

$$Z_n = b_n C_\alpha^{1/\alpha} \max_{k=0,\ldots,n-1} \left| \Gamma_1^{-1/\alpha} \frac{|f_k(U_1^{(n)})|}{\max_{i=0,\ldots,n-1} |f_i(U_1^{(n)})|} \right|$$

and $K_0$ be the smallest $k = 0, 1, \ldots, n-1$ over which the maximum is achieved, we see that both $Z_n$ and $K_0$ are measurable $\mathcal{G}$. Further, the symmetry tells us that, for any $x > 0$,

$$P(X_{K_0} > x | \mathcal{G}) \geq \tfrac{1}{2} P(Z_n > x | \mathcal{G}).$$

Hence, for any $x > 0$,

$$P\left(\frac{1}{c_n} M_n > x\right)$$

(4.13)
$$\geq \frac{1}{2} P\left( \max_{k=0,\ldots,n-1} b_n C_\alpha^{1/\alpha} \left| \varepsilon_1 \Gamma_1^{-1/\alpha} \frac{|f_k(U_1^{(n)})|}{\max_{i=0,\ldots,n-1} |f_i(U_1^{(n)})|} \right| > c_n x \right)$$
$$= \frac{1}{2} P\left( \Gamma_1^{-1/\alpha} > C_\alpha^{-1/\alpha} \frac{c_n}{b_n} x \right) \to \frac{1}{2}$$

as $n \to \infty$. Hence lack of tightness.

Suppose now that (4.6) holds. Let $K = 1, 2, \ldots$ be such that

(4.14)
$$\alpha (K+1) \theta > 1.$$

We claim that, in this case, for all $\epsilon > 0$ satisfying

(4.15)
$$0 < \epsilon < \frac{1}{K},$$

we have

(4.16)
$$P\left( \max_{k=0,\ldots,n-1} |X_k| > \lambda b_n, \Gamma_1^{-1/\alpha} \leq \epsilon \lambda \right) \to 0 \quad \text{as } n \to \infty$$

for all $\lambda > 0$. Indeed, choose

(4.17)
$$\frac{1}{\theta} < p < \alpha(K+1).$$

Notice that the probability in the left-hand side of (4.16) is bounded from above by

$$\sum_{k=0}^{n-1} P\left( |X_k| > \lambda b_n, \Gamma_j^{-1/\alpha} \frac{|f_k(U_j^{(n)})|}{\max_{i=0,\ldots,n-1} |f_i(U_j^{(n)})|} \leq \epsilon \lambda \text{ for all } j = 1, 2, \ldots \right).$$

For $f$ in (2.2), let

$$\|f\|_\alpha = \left(\int_E |f(x)|^\alpha \, m(dx)\right)^{1/\alpha},$$

and notice that, for any $k = 0, 1, \ldots, n-1$, the points

$$b_n \varepsilon_j \Gamma_j^{-1/\alpha} \frac{f_k(U_j^{(n)})}{\max_{i=0,\ldots,n-1} |f_i(U_j^{(n)})|}, \qquad j = 1, 2, \ldots,$$

represent a symmetric Poisson random measure on $\mathbb{R}$ whose mean measure assigns a mass of $x^{-\alpha}\|f\|_\alpha^\alpha/2$ to the set $(x, \infty)$ for every $x > 0$ [see, e.g., Propositions 4.3.1 and 4.4.1 in Resnick (1992)]. Since the same random measure can be represented by the points

$$\varepsilon_j \Gamma_j^{-1/\alpha} \|f\|_\alpha, \qquad j = 1, 2, \ldots,$$

we conclude that the probability in (4.16) is bounded from above by

$$nP\left(C_\alpha^{1/\alpha} \sum_{j=1}^\infty \varepsilon_j \Gamma_j^{-1/\alpha} > \lambda\|f\|_\alpha^{-1} b_n, \Gamma_j^{-1/\alpha} \leq \epsilon\lambda\|f\|_\alpha^{-1} b_n\right)$$

$$\leq nP\left(C_\alpha^{1/\alpha} \sum_{j=K+1}^\infty \varepsilon_j \Gamma_j^{-1/\alpha} > \lambda(1-\epsilon K)\|f\|_\alpha^{-1} b_n\right)$$

$$\leq n b_n^{-p} \frac{\|f\|_\alpha^p E|C_\alpha^{1/\alpha} \sum_{j=K+1}^\infty \varepsilon_j \Gamma_j^{-1/\alpha}|^p}{\lambda^p(1-\epsilon K)^p}.$$

As long as the expectation above is finite, the latter expression goes to 0 as $n \to \infty$ and, hence, (4.16) follows. The expectation is finite by the choice of $p$ in (4.17). Indeed, notice that $E\Gamma_j^{-p/\alpha} < \infty$ for all $j \geq K+1$ and that, by the Stirling formula, $E\Gamma_j^{-p/\alpha} \sim e^{p/\alpha} j^{-p/\alpha}$ as $j \to \infty$. Assuming without loss of generality that $p/2 = m$ is an integer (we can remove finitely many leading terms in the sum and increase $p$, if necessary), we see that for finite positive constants $c_1, c_2$,

$$E\left|\sum_{j=K+1}^\infty \varepsilon_j \Gamma_j^{-1/\alpha}\right|^p \leq c_1 E\left(\sum_{j=K+1}^\infty \Gamma_j^{-2/\alpha}\right)^{p/2}$$

$$= c_1 \sum_{j_1=K+1}^\infty \cdots \sum_{j_m=K+1}^\infty E \prod_{i=1}^m \Gamma_{j_i}^{-2/\alpha}$$

$$\leq c_1 \left(\sum_{j=K+1}^\infty (E\Gamma_j^{-2m/\alpha})^{1/m}\right)^m$$

EXTREMES OF STATIONARY STABLE PROCESSES    15

$$= c_1 \left( \sum_{j=K+1}^{\infty} (E\Gamma_j^{-p/\alpha})^{2/p} \right)^{p/2}$$

$$\leq c_2 \left( \sum_{j=K+1}^{\infty} j^{-2/\alpha} \right)^{p/2} < \infty.$$

Fix $\epsilon > 0$ satisfying (4.15). Given $\delta > 0$, choose $\lambda > 0$ such that $P(\Gamma_1^{-1/\alpha} > \epsilon\lambda) \leq \delta/2$, choose $n_0$ such that

$$P\left( \max_{k=0,\ldots,n-1} |X_k| > \lambda b_n, \Gamma_1^{-1/\alpha} \leq \epsilon\lambda \right) \leq \frac{\delta}{2} \qquad \text{for } n > n_0$$

and $\lambda' \geq \lambda$ such that

$$P(M_k > \lambda' b_k) \leq \delta \qquad \text{for } k = 1, \ldots, n_0.$$

Then

$$P\left( \frac{1}{b_n} M_n > \lambda' \right) \leq \delta$$

for all $n \geq 1$, and so (4.5) holds with $f_n \equiv 1$.

Now, suppose that **X** is generated by a conservative flow. Let **Y** be a stationary S$\alpha$S process independent of **X**, also given by an integral representation of the form (2.1), say,

$$Y_n = \int_{E'} g_n(x) M'(dx), \qquad n = 0, 1, 2, \ldots,$$

where $M'$ is a S$\alpha$S random measure with control measure $m'$, independent of $M$ in the integral representation of **X**, with the functions $g_n$ also given in the form (2.2), with some nonsingular conservative flow $\phi'$ on $E'$, and such that

$$b_n^Y = \left( \int_{E'} \max_{j=0,\ldots,n-1} |g_j(x)|^\alpha m'(dx) \right)^{1/\alpha}, \qquad n = 1, 2, \ldots,$$

satisfies (4.6) for some $\theta > 0$. Processes **Y** with the above properties exist; see the examples in the next section. However, the above step may require enlarging the probability space we are working with. Let $\mathbf{Z} = \mathbf{X} + \mathbf{Y}$. Then **Z** is a stationary S$\alpha$S process generated by a conservative flow. We use its natural integral representation on $E \cup E'$ with the naturally defined flow on that space. Let $b_n^Z$ be the corresponding quantity in (3.1) defined for the process **Z**. Note that $b_n^Z \geq b_n^Y$ for all $n$, hence the process **Z** satisfies (4.6) as well. By the already proven part of the theorem, the sequence $(b_n^Z)^{-1} \max_{k=0,\ldots,n-1} |Z_k|$, $n = 1, 2, \ldots$, is tight. Since, for any $x > 0$ and $n = 1, 2, \ldots$,

$$P\left( \max_{k=0,\ldots,n-1} |Z_k| > x \right) \geq \tfrac{1}{2} P\left( \max_{k=0,\ldots,n-1} |X_k| > x \right)$$

by the symmetry of **Y**, we conclude that the sequence $(b_n^Z)^{-1} \max_{k=0,\ldots,n-1} |X_k|$, $n=1,2,\ldots$, is tight as well.

However, the process **Z** is generated by a conservative flow and, hence, by Theorem 3.1, $b_n^Z = o(n^{1/\alpha})$. Therefore, (4.3) follows.

Suppose now that (4.8) holds. Then for every $1 \leq j_1 < j_2$ and $\epsilon > 0$

$$P\left(\text{for some } k=0,1,\ldots,n-1,\ \Gamma_{j_i}^{-1/\alpha} \frac{|f_k(U_{j_i}^{(n)})|}{\max_{m=0,\ldots,n-1} |f_m(U_{j_i}^{(n)})|} > \epsilon,\ i=1,2\right)$$

$$\leq P(\Gamma_1 \leq \tau) + P\left(\text{for some } k=0,1,\ldots,n-1,\right.$$

$$\left.\frac{|f_k(U_j^{(n)})|}{\max_{m=0,\ldots,n-1} |f_m(U_j^{(n)})|} > \epsilon \tau^{1/\alpha},\ j=1,2\right)$$

for any $\tau > 0$. Letting first $n \to \infty$ and then $\tau \to 0$ shows that, for every $1 \leq j_1 < j_2$ and $\epsilon > 0$,

(4.18)
$$\lim_{n\to\infty} P\left(\text{for some } k=0,1,\ldots,n-1,\right.$$
$$\left.\Gamma_{j_i}^{-1/\alpha} \frac{|f_k(U_{j_i}^{(n)})|}{\max_{m=0,\ldots,n-1} |f_m(U_{j_i}^{(n)})|} > \epsilon,\ i=1,2\right) = 0.$$

Observe, further, that for any $\epsilon > 0$,

$$P\left(\text{for some } k=0,1,\ldots,n-1,\right.$$

$$\left.\Gamma_j^{-1/\alpha} \frac{|f_k(U_j^{(n)})|}{\max_{m=0,\ldots,n-1} |f_m(U_j^{(n)})|} > \epsilon \text{ for at least 2 different } j\right)$$

$$=: \varphi_n^{(1)}(\epsilon) \leq P(\Gamma_J^{-1/\alpha} > \epsilon)$$

$$+ \sum_{j_1=1}^{J-1} \sum_{j_2=j_1+1}^{J-1} P\left(\text{for some } k=0,1,\ldots,n-1,\right.$$

$$\left.\Gamma_{j_i}^{-1/\alpha} \frac{|f_k(U_{j_i}^{(n)})|}{\max_{m=0,\ldots,n-1} |f_m(U_{j_i}^{(n)})|} > \epsilon,\ i=1,2\right)$$

for any $J = 1, 2, \ldots$. Letting $n \to \infty$ and using (4.18), and then letting $J \to \infty$ shows that, for every $\epsilon > 0$

(4.19) $$\lim_{n\to\infty} \varphi_n^{(1)}(\epsilon) = 0.$$

Suppose now that both (4.6) and (4.8) hold. Let $K$ be as in (4.14). Let $\epsilon > 0$ and $0 < \delta < 1$ satisfy

(4.20) $$0 < \epsilon < \frac{\delta}{K}.$$

For any $\lambda > 0$, we have

$$P\left(\frac{1}{b_n} M_n > \lambda\right)$$
$$\leq P(C_\alpha^{1/\alpha} \Gamma_1^{-1/\alpha} > \lambda(1-\delta)) + \varphi_n^{(1)}(C_\alpha^{-1/\alpha} \epsilon \lambda)$$
$$+ P\left(\max_{k=0,\ldots,n-1} \left|\sum_{j=1}^\infty \varepsilon_j \Gamma_j^{-1/\alpha} \frac{f_k(U_j^{(n)})}{\max_{i=0,\ldots,n-1} |f_i(U_j^{(n)})|}\right| > C_\alpha^{-1/\alpha} \lambda,\right.$$

(4.21) $$\Gamma_1^{-1/\alpha} \leq C_\alpha^{-1/\alpha} \lambda(1-\delta),$$

and for each $k = 0, 1, \ldots, n-1$,

$$\Gamma_j^{-1/\alpha} \frac{|f_k(U_j^{(n)})|}{\max_{i=0,\ldots,n-1} |f_i(U_j^{(n)})|} > C_\alpha^{-1/\alpha} \epsilon \lambda$$

$$\left. \text{for at most one } j = 1, 2, \ldots \right)$$

$$=: P(C_\alpha^{1/\alpha} \Gamma_1^{-1/\alpha} > \lambda(1-\delta)) + \varphi_n^{(1)}(C_\alpha^{-1/\alpha} \epsilon \lambda) + \varphi_n^{(2)}(\epsilon, \delta).$$

Proceeding similarly to the argument used in proving (4.16) we have

$$\varphi_n^{(2)}(\epsilon, \delta) \leq \sum_{k=0}^{n-1} P\left(\left|\sum_{j=1}^\infty \varepsilon_j \Gamma_j^{-1/\alpha} \frac{f_k(U_j^{(n)})}{\max_{i=0,\ldots,n-1} |f_i(U_j^{(n)})|}\right| > C_\alpha^{-1/\alpha} \lambda,\right.$$

$$\Gamma_j^{-1/\alpha} \frac{|f_k(U_j^{(n)})|}{\max_{i=0,\ldots,n-1} |f_i(U_j^{(n)})|}$$
$$\leq C_\alpha^{-1/\alpha} \lambda(1-\delta) \text{ for each } j = 1, 2, \ldots$$

$$\text{and } \Gamma_j^{-1/\alpha} \frac{|f_k(U_j^{(n)})|}{\max_{i=0,\ldots,n-1} |f_i(U_j^{(n)})|} > C_\alpha^{-1/\alpha} \epsilon \lambda$$

(4.22)
$$\left. \text{for at most one } j = 1, 2, \ldots \right)$$

$$\leq nP\left(\left|\sum_{j=1}^\infty \varepsilon_j \Gamma_j^{-1/\alpha}\right| > C_\alpha^{-1/\alpha} \lambda \|f\|_\alpha^{-1} b_n,\right.$$

$$\Gamma_1^{-1/\alpha} \leq C_\alpha^{-1/\alpha}\lambda(1-\delta)\,\|f\|_\alpha^{-1}b_n,$$

$$\text{and } \Gamma_2^{-1/\alpha} \leq C_\alpha^{-1/\alpha}\epsilon\lambda\|f\|_\alpha^{-1}b_n\Bigg)$$

and the latter expression goes to zero as $n \to \infty$ by the choice of $\epsilon$ and $\delta$, as in the proof of (4.16).

We conclude by (4.19) and (4.22) that, for any $0 < \delta < 1$,

$$\limsup_{n\to\infty} P\Big(\frac{1}{b_n}M_n > \lambda\Big) \leq P(C_\alpha^{1/\alpha}\Gamma_1^{-1/\alpha} > \lambda(1-\delta))$$

$$= 1 - \exp\{-C_\alpha\lambda^{-\alpha}(1-\delta)^{-\alpha}\},$$

and letting $\delta \to 0$, we obtain

(4.23) $$\limsup_{n\to\infty} P\Big(\frac{1}{b_n}M_n > \lambda\Big) \leq 1 - \exp\{-C_\alpha\lambda^{-\alpha}\}.$$

In the opposite direction, the argument is similar. For any $\epsilon$ and $\delta > 0$ satisfying (4.20), we have

$$P\Big(\frac{1}{b_n}M_n > \lambda\Big)$$

$$\geq P(C_\alpha^{1/\alpha}\Gamma_1^{-1/\alpha} > \lambda(1+\delta)) - \varphi_n^{(1)}(C_\alpha^{-1/\alpha}\epsilon\lambda)$$

$$- P\Bigg(\max_{k=0,\ldots,n-1}\bigg|\sum_{j=1}^\infty \varepsilon_j\Gamma_j^{-1/\alpha}\frac{f_k(U_j^{(n)})}{\max_{i=0,\ldots,n-1}|f_i(U_j^{(n)})|}\bigg| \leq C_\alpha^{-1/\alpha}\lambda,$$

$$\Gamma_1^{-1/\alpha} > C_\alpha^{-1/\alpha}\lambda(1+\delta),$$

and for each $k = 0, 1, \ldots, n-1$,

$$\Gamma_j^{-1/\alpha}\frac{|f_k(U_j^{(n)})|}{\max_{i=0,\ldots,n-1}|f_i(U_j^{(n)})|} > C_\alpha^{-1/\alpha}\epsilon\lambda$$

$$\text{for at most one } j = 1, 2, \ldots\Bigg)$$

$$=: P(C_\alpha^{1/\alpha}\Gamma_1^{-1/\alpha} > \lambda(1+\delta)) - \varphi_n^{(1)}(C_\alpha^{-1/\alpha}\epsilon\lambda) - \varphi_n^{(3)}(\epsilon,\delta).$$

Once again, the choice of $\epsilon$ and $\delta$ gives us

(4.24) $$\lim_{n\to\infty} \varphi_n^{(3)}(\epsilon,\delta) = 0,$$

and so we conclude by (4.19) and (4.24), that for any $\delta > 0$,

$$\liminf_{n\to\infty} P\Big(\frac{1}{b_n}M_n > \lambda\Big) \geq P(C_\alpha^{1/\alpha}\Gamma_1^{-1/\alpha} > \lambda(1+\delta))$$

$$= 1 - \exp\{-C_\alpha\lambda^{-\alpha}(1+\delta)^{-\alpha}\},$$

and letting $\delta \to 0$, we obtain a lower bound matching (4.23). This proves (4.9).

If **X** is not generated by a conservative flow, then it follows from Theorem 3.1 that
$$n^{-1/\alpha} b_n \to K_X \qquad \text{as } n \to \infty.$$
In particular, both conditions (4.6) and (4.8) are satisfied (see Remark 4.2). Therefore, (4.1) follows from the already proven (4.9), and the proof of all parts of the theorem is complete. □

REMARK 4.3. The entire statement of part (ii) of Theorem 4.1 remains valid with any integral representation of the process, not necessarily with functions of the form (2.2). Of course, the sequence $(b_n)$ does not depend on the representation. One can easily see that the particular structure of the integral representation was not used in the proof, except when proving (4.3). The latter statement is, however, a distributional one, and does not depend on the integral representation.

REMARK 4.4. The assumptions (4.6) and (4.8) mean, intuitively, that one and only one Poisson jump in the series representation (4.12) significantly contributes to the value of $M_n$ for large $n$. Because of that, an extreme value distribution arises as a limit. Either of these two assumptions may fail, as will be seen from the examples in the next section. Even though a complete limit theory in such cases is unavailable at the moment, limiting distributions (when weak limits exist) are likely to depend on the number of Poisson jumps that contribute significantly to the value of the maximum. In particular, the limiting distribution is not, in general, an extreme value distribution. See Example 5.1.

In any case, no subsequential weak limit of $M_n/b_n$ can be constant, as can be seen by using (4.13) with $c_n = b_n$:
$$P\left(\frac{1}{b_n} M_n > x\right) \geq \frac{1}{2} P(\Gamma_1^{-1/\alpha} > C_\alpha^{-1/\alpha} x)$$
for all $n \geq 1$, so that any subsequential weak limit must have a nonvanishing tail.

REMARK 4.5. Given a stationary S$\alpha$S, process **X**, let

(4.25) $\qquad M_n^{(0)} = \max(X_0, X_1, \ldots, X_{n-1}), \qquad n = 1, 2, \ldots.$

The same argument as that used in proof of Theorem 4.1 shows that, if **X** is not generated by a conservative flow, then

(4.26) $\qquad \dfrac{1}{n^{1/\alpha}} M_n^{(0)} \Rightarrow C_\alpha^{1/\alpha} K_X^{(0)} Z_\alpha$

weakly as $n \to \infty$, where

$$K_X^{(0)} = \left( \tfrac{1}{2} \int_W g_+^{(0)}(v)^\alpha \, \nu(dv) + \tfrac{1}{2} \int_W g_-^{(0)}(v)^\alpha \, \nu(dv) \right)^{1/\alpha},$$

with

$$g_+^{(0)}(v) = \sup_{k=0,\pm 1, \pm 2,\ldots} f(v,k)_+, \qquad g_-^{(0)}(v) = \sup_{k=0,\pm 1, \pm 2,\ldots} f(v,k)_-,$$

$v \in W$, for any representation of $\mathbf{X}^C$ in the mixed moving average form (2.4). Here $a_+$ and $a_-$ are, correspondingly, the positive and negative parts of a real number $a$. In the particular case of linear $\alpha$-stable process, this result is in Theorem 3.8.3 in Leadbetter, Lindgren and Rootzén (1983).

It is, of course, clear that if $\mathbf{X}$ is generated by a conservative flow, then $n^{-1/\alpha} M_n^{(0)} \to 0$ in probability as $n \to \infty$.

Here is a sketch of an argument for (4.26). The points

$$\left( C_\alpha^{1/\alpha} b_n \, \varepsilon_j \, \Gamma_j^{-1/\alpha} \frac{\mathbf{f}(U_j^{(n)})}{\max_{i=0,\ldots,n-1} |f_i(U_j^{(n)})|} \right), \qquad j = 1, 2, \ldots,$$

in (4.12) form a symmetric Poisson random measure on $\mathbb{R}^n$. Here $\mathbf{f}(x) = (f_0(x), \ldots, f_{n-1}(x))$. The proof of Theorem 4.1 shows that the event $\{n^{-1/\alpha} \times M_n^{(0)} > \lambda\}$ is asymptotically equivalent to the event that at least one of the points of the Poisson random measure is in the set $((-\infty, \lambda n^{1/\alpha}]^n)^c$.

It is easy to check that the mean measure of that set is equal to

$$C_\alpha \lambda^{-\alpha} \frac{1}{n} \left( \frac{1}{2} b_{n,+}^\alpha + \frac{1}{2} b_{n,-}^\alpha \right),$$

where

$$b_{n,\pm} = \left( \int_E \max_{j=0,\ldots,n-1} (f_j(x))_\pm^\alpha m(dx) \right)^{1/\alpha}, \qquad n = 1, 2, \ldots.$$

An argument identical to that used to prove the second part of Theorem 3.1 shows that

$$\frac{1}{n} \left( \frac{1}{2} b_{n,+}^\alpha + \frac{1}{2} b_{n,-}^\alpha \right) \to (K_X^{(0)})^\alpha \qquad \text{as } n \to \infty$$

and (4.26) follows.

REMARK 4.6. Theorem 4.1 extends easily to the case of complex-valued rotationally invariant stationary stable processes; we refer the reader to Chapter 6 in Samorodnitsky and Taqqu (1994) for basic information on such processes. Such processes also have an integral representation of the form (2.1) and (2.2), but this time $M$ is a complex-valued rotationally invariant S$\alpha$S random measure, and the function $f$ in (2.2) is complex-valued as well. Further results, like (2.3) and (2.4), hold in the complex-valued case as well. See Rosiński (1995).

In the complex-valued case, one has to replace the series representation (4.12) by

$$X_k = b_n\, C_\alpha^{1/\alpha} \sum_{j=1}^{\infty} e^{i A_j}\, \Gamma_j^{-1/\alpha}\, \frac{f_k(U_j^{(n)})}{\max_{i=0,\ldots,n-1}|f_i(U_j^{(n)})|}, \quad (4.27)$$
$$k = 0, 1, \ldots, n-1,$$

where $(\Gamma_j)$ and $(U_j^{(n)})$ are as before, and $(A_j)$ is an independent of them i.i.d. sequence of random variables uniformly distributed in $(0, 2\pi)$. This representation can be easily derived from the real-valued series representations in Chapter 3 of Samorodnitsky and Taqqu (1994).

It can be easily verified that the representation (4.27) is a perfect substitute for (4.12) in the proof of Theorem 4.1 and, hence, all the claims of that theorem hold in the rotationally invariant complex-valued case as well.

REMARK 4.7. An open, and very interesting, question is what are some possible counterparts of Theorem 4.1 in the case of continuous-time $\alpha$-stable processes.

Integral representations of the type (2.1) with the functions $(f_n)$ given in the form (2.2) exist not only for S$\alpha$S processes, but also for all strictly $\alpha$-stable processes with $\alpha \neq 1$ (in this case the random measure $M$ may not be symmetric, and will, in general, have a nonconstant skewness intensity). See Rosiński (1994). The following result is a straightforward consequence of Theorem 4.1.

THEOREM 4.8. *Let* $\mathbf{X} = (X_0, X_1, \ldots)$ *be a stationary strictly $\alpha$-stable process, $0 < \alpha < 2$, $\alpha \neq 1$, with integral representation* (2.1), *and functions* $(f_n)$ *given by* (2.2).

(i) *Suppose that* $\mathbf{X}$ *is not generated by a conservative flow. Then* (4.1) *holds.*

(ii) *Suppose that* $\mathbf{X}$ *is generated by a conservative flow. Then* (4.3) *holds.*

PROOF. In the case $0 < \alpha < 1$, a strictly $\alpha$-stable process has a series representation of the form (4.12) [but with $(\varepsilon_j)$ dependent on $(U_j^{(n)})$ and not necessarily symmetric]. The entire argument of Theorem 4.1 used to prove (4.1) and (4.3) goes through in this case. Hence, we only need to consider the case $1 < \alpha < 2$. In that case there is a positive constant $\tau > 0$ such that

$$\min(P(X_1 > 0), P(X_1 < 0)) \geq \tau. \quad (4.28)$$

See Theorem 4 in Zolotarev (1957).

Let **Y** be an independent copy of **X**. Then $\mathbf{Z} = \mathbf{X} - \mathbf{Y}$ is a stationary S$\alpha$S process that has an integral representation (2.1) with the same functions $(f_n)$ as **X** does, but where the random measure $M$ is now replaced by a S$\alpha$S random measure $M'$ whose control measure is given by $m' = 2m$; see Samorodnitsky and Taqqu (1994). If, $(b_n^Z)$, is the sequence (3.1) defined for **Z** using this integral representation, we have

(4.29) $$b_n^Z = 2^{1/\alpha} b_n, \qquad n = 1, 2, \ldots.$$

Suppose that **X** is generated by a conservative flow; then so is **Z**. We have by (4.28), for any $\lambda > 0$,

$$P\left(\max_{k=0,\ldots,n-1} |Z_n| > \lambda\right) \geq P(|Z_{K_n}| > \lambda)$$
$$\geq P(X_{K_n} > \lambda, Y_{K_n} > 0) + P(X_{K_n} < -\lambda, Y_{K_n} < 0)$$
$$\geq \tau P(|X_{K_n}| > \lambda) = \tau P\left(\max_{k=0,\ldots,n-1} |X_n| > \lambda\right),$$

by the independence of **X** and **Y**, where $K_n$ is the smallest $k = 0, 1, \ldots$ such that $|X_k| = \max_{j=0,\ldots,n-1} |X_j|$. Therefore, (4.3) for **X** follows from (4.3) for **Z** proven in Theorem 4.1.

If **X** is not generated by a conservative flow, then neither is **Z**. We have by Theorem 4.1 and (4.29),

(4.30) $$\frac{1}{n^{1/\alpha}} \max_{k=0,\ldots,n-1} |Z_n| \Rightarrow 2^{1/\alpha} C_\alpha^{1/\alpha} K_X Z_\alpha$$

weakly as $n \to \infty$.

Let

$$K_X(n) = \min\{k = 0, 1, \ldots, n-1 : |X_k| = M_n\}$$

and similarly with $K_Y(n)$. For every $\lambda > 0$ and $\epsilon > 0$, we have, by the stationarity and independence,

$$P\left(\max_{k=0,\ldots,n-1} |Z_n| \leq \lambda n^{1/\alpha}\right)$$
$$= P\left(\max_{k=0,\ldots,n-1} |X_n - Y_n| \leq \lambda n^{1/\alpha}\right)$$
$$\leq P\left(\max_{k=0,\ldots,n-1} |X_n| \leq \lambda(1+\epsilon)n^{1/\alpha}\right) P\left(\max_{k=0,\ldots,n-1} |Y_n| \leq \lambda(1+\epsilon)n^{1/\alpha}\right)$$
$$\quad + P(|Y_{K_X(n)}| > \epsilon n^{1/\alpha}, \text{ or } |X_{K_Y(n)}| > \epsilon n^{1/\alpha})$$
$$\leq \left(P\left(\max_{k=0,\ldots,n-1} |X_n| \leq \lambda(1+\epsilon)n^{1/\alpha}\right)\right)^2 + 2P(|X_1| > \epsilon n^{1/\alpha}).$$

Therefore, using (4.30), we have, for any $\lambda > 0$,

$$\liminf_{n\to\infty} P\left(\max_{k=0,\ldots,n-1} |X_n| \leq \lambda n^{1/\alpha}\right) \geq \exp\{-C_\alpha \lambda^{-\alpha}(1+\epsilon)^\alpha\}$$

for any $\epsilon > 0$. Letting $\epsilon \to 0$, we obtain

(4.31) $$\liminf_{n\to\infty} P\left(\max_{k=0,\ldots,n-1} |X_n| \leq \lambda n^{1/\alpha}\right) \geq \exp\{-C_\alpha \lambda^{-\alpha}\}$$

for any $\lambda > 0$.

In the other direction, for every $\lambda > 0$ and $0 < \epsilon < 1$

$$P\left(\max_{k=0,\ldots,n-1} |Z_n| \leq \lambda n^{1/\alpha}\right)$$
$$\geq \left(P\left(\max_{k=0,\ldots,n-1} |X_n| \leq \lambda(1-\epsilon)n^{1/\alpha}\right)\right)^2$$
$$- P(\text{for some } k = 0, 1, \ldots, n-1, |X_k| > \epsilon n^{1/\alpha} \text{ and } |Y_k| > \epsilon n^{1/\alpha}).$$

By stationarity and independence,

$$P(\text{for some } k = 0, 1, \ldots, n-1, |X_k| > \epsilon n^{1/\alpha} \text{ and } |Y_k| > \epsilon n^{1/\alpha})$$
$$\leq n(P(|X_1| > \epsilon n^{1/\alpha}))^2 \to 0$$

as $n \to \infty$. Once again, using (4.30), we have, for any $\lambda > 0$,

$$\limsup_{n\to\infty} P\left(\max_{k=0,\ldots,n-1} |X_n| \leq \lambda n^{1/\alpha}\right) \leq \exp\{-C_\alpha \lambda^{-\alpha}(1-\epsilon)^\alpha\}$$

for any $0 < \epsilon < 1$. Letting $\epsilon \to 0$, we obtain an upper bound matching (4.31) and, hence, complete the proof. $\square$

REMARK 4.9. We see immediately from the proof of Theorem 4.8 that the statement (4.9) extends to the skewed case as well if (4.10) holds.

**5. Examples.** The results of the previous section describe completely the limiting behavior of the partial maxima of stationary S$\alpha$S processes not generated by a conservative flow. The picture for processes generated by conservative flows is less complete, and in this section we consider several examples of such processes to illustrate what may happen.

Our first example shows that, in general, the partial maximum $M_n$ does not have an extreme value limit.

EXAMPLE 5.1. Let $Z_0, Z_1, \ldots$ be i.i.d. standard normal random variables, independent of a positive $(\alpha/2)$-stable random variable $A$ with Laplace transform $Ee^{-\theta A} = e^{-\theta^{\alpha/2}}, \theta \geq 0$. Then $X_n = A^{1/2} Z_n, n = 0, 1, 2, \ldots$, is a stationary S$\alpha$S process, the simplest type of *sub-Gaussian* S$\alpha$S processes; see Section 3.7 in Samorodnitsky and Taqqu (1994). This process has an integral representation of the form

$$(5.1) \qquad X_n = (d_\alpha)^{-1} \int_{\mathbb{R}^{\mathbb{Z}}} g_n \, M(d\mathbf{g}), \qquad n = 0, 1, 2, \ldots,$$

where $d_\alpha = \sqrt{2}(E|Z_0|^\alpha)^{1/\alpha}$, and $M$ is an S$\alpha$S random measure on $\mathbb{R}^{\mathbb{Z}}$ whose control measure $m$ is a probability measure under which the projections $g_k$, $k \in \mathbb{Z}$, are i.i.d. standard normal random variables. With $\phi$ being the usual left shift operator on $\mathbb{R}^{\mathbb{Z}}$, this is a representation with the functions given in the form (2.2), with $a_n \equiv 1$ and a measure-preserving flow. The flow is, obviously, conservative, as is any measure-preserving flow on a finite measure space. As usual, the representation (5.1) is in law.

An elementary direct computation shows that

$$b_n^\alpha = (d_\alpha)^{-\alpha} E \max_{k=0,1,\ldots,n-1} |Z_k|^\alpha \sim (d_\alpha)^{-\alpha}(2 \log n)^{\alpha/2}$$

and

$$\frac{1}{(2 \log n)^{1/2}} \max_{k=0,1,\ldots,n-1} |Z_k| \to 1 \qquad \text{with probability 1}$$

as $n \to \infty$. Therefore, the assumption (4.6) in Theorem 4.1 fails, and we see directly that

$$\frac{1}{b_n} M_n \Rightarrow (d_\alpha)^{-1} A^{1/2},$$

a nonextreme value limit, even though the partial maximum $M_n$ does grow at the rate $b_n$.

Our next example shows that, if $1 \leq \alpha < 2$ and without the assumption (4.6), the sequence of the partial maxima of the process may grow faster than the sequence $(b_n)$; that is, (4.5) may not hold with $f_n \equiv 1$.

EXAMPLE 5.2. Let **X** be given by (5.1) with $d_\alpha = 1$, where this time the control measure $m$ of the S$\alpha$S random measure $M$ is a probability measure under which the projections $g_k$, $k = \ldots, -1, 0, 1, 2, \ldots$ are i.i.d. Rademacher random variables. We have, once again, a conservative measure-preserving flow, and $b_n \equiv 1$. Let $1 \leq \alpha < 2$. We claim that (4.5) with $f_n \equiv 1$ does not hold. Indeed, if it did, then the S$\alpha$S process **X** would be a.s. bounded, and then, for any $0 < p < \alpha$, we would have

$$EM_\infty^p =: E \sup_{k=0,1,\ldots} |X_k|^p < \infty$$

[see, e.g., Araujo and Giné (1980)], which would contradict, in the case $1 < \alpha < 2$, (2.22) in Marcus (1984) and in the case $\alpha = 1$, it would contradict (2.23) in that paper.

The next example exhibits a variety of rates at which the sequence $(b_n)$ and the partial maxima of the processes can grow. We look at a class of stationary S$\alpha$S processes generated by null recurrent Markov chains, introduced in Rosiński and Samorodnitsky (1996) and studied in more details in Resnick, Samorodnitsky and Xue (2000) and Mikosch and Samorodnitsky (2000).

EXAMPLE 5.3. We start with an irreducible null-recurrent Markov chain on $\mathbb{Z}$ with law $P_i(\cdot)$, $i \in \mathbb{Z}$ and transition probabilities $(p_{jk})$. Let $\pi = (\pi_i)_{i \in \mathbb{Z}}$ be the $\sigma$-finite invariant measure corresponding to the family $(P_i)$ satisfying $\pi_0 = 1$. Let $P_i^*$ be the bilateral extension of $P_i$ to $E = \mathbb{Z}^{\mathbb{Z}}$; that is, under $P_i^*$, $x_0 = i$, $(x_0, x_1, \ldots)$ is a Markov chain with transition probabilities $(p_{jk})$ and $(x_0, x_{-1}, \ldots)$ is a Markov chain with transition probabilities $(\pi_k p_{kj}/\pi_j)$. Define a $\sigma$-finite measure $m$ on $E$ by

$$m(\cdot) = \sum_{i=-\infty}^{\infty} \pi_i P_i^*(\cdot),$$

and observe that $m$ is invariant under the left shift operator $\phi$, and the latter generates a conservative flow [see Harris and Robbins (1953)].

Let $\mathbf{X}$ be a stationary S$\alpha$S process defined by the integral representation (2.1), with $M$ being an S$\alpha$S random measure with control measure $m$,

$$f_n(x) = f \circ \phi^n(x), \qquad x \in E, \ n = 0, 1, 2, \ldots,$$

with

$$f(x) = \mathbf{1}_0(x_0), \qquad x = (\ldots, x_{-1}, x_0, x_1, x_2, \ldots) \in E.$$

Because of the null recurrence of the Markov chain, this is a mixing stationary process (unlike, say, the processes in Examples 5.1 and 5.2); see Rosiński and Samorodnitsky (1996).

Let $d$ be the period of the Markov chain [the largest common factor of $n \geq 1$ such that $P_0(x_n = 0) > 0$]. Assume that

(5.2) $$P_0(x_{nd} = 0) = n^{-\gamma} L(n) \qquad \text{as } n \to \infty$$

for some $\gamma \in (0, 1)$ and a slowly varying function $L$. Let

$$\tau = \tau(x) = \inf\{n \geq 1 : x_n = 0\}$$

be the first return time to 0. Notice that

$$b_n^\alpha = \sum_{i=-\infty}^{\infty} \pi_i P_i(x_k = 0 \text{ for some } k = 0, \ldots, n-1)$$

$$= P_0(\tau \geq n) + \sum_{i=-\infty}^{\infty} \pi_i P_i(\tau \leq n-1)$$

$$= P_0(\tau \geq n) + m(\tau \leq n-1)$$

$$\sim \frac{d^{1-\gamma}}{\Gamma(1+\gamma)\Gamma(1-\gamma)} n^\gamma (L_n)^{-1}$$

as $n \to \infty$ by Remark 3.1 and Lemma 3.3 in Resnick, Samorodnitsky and Xue (2000). In particular, it follows from Theorem 4.1 and Remark 4.2 that

$$(5.3) \qquad \frac{L(n)^{1/\alpha}}{n^{\gamma/\alpha}} M_n \Rightarrow \left(\frac{C_\alpha d^{1-\gamma}}{\Gamma(1+\gamma)\Gamma(1-\gamma)}\right)^{1/\alpha} Z_\alpha$$

as $n \to \infty$ if $\gamma > 1/2$.

What happens if $\gamma \leq 1/2$? Let $Y_0^{(i)}, Y_1^{(i)}, \ldots, i = 1, 2$, be two independent Markov chains with the same transition probabilities $(p_{jk})$ as before. Then

$$Y_n^* = (Y_n^{(1)}, Y_n^{(2)}), \qquad n = 0, 1, 2 \ldots,$$

is a Markov chain with state space $\mathbb{Z}^2$ and transition probabilities

$$p^*_{(j_1,j_2),(k_1,k_2)} = p_{j_1,k_1} p_{j_2,k_2}, \qquad j_1, j_2, k_1, k_2 \in \mathbb{Z}.$$

Let

$$\tau^* = \inf\{n \geq 1 : Y_n^* = (0,0)\}$$

be the first time the new Markov chain returns to $(0,0)$.

Let $r_n$ be the probability in the left-hand side of (4.8). For any $0 < \epsilon < 1$, we have

$$r_n = b_n^{-2\alpha}\left[P_{(0,0)}(\tau^* \geq n) + \sum_{i=-\infty}^{\infty}\sum_{j=-\infty}^{\infty} \pi_i \pi_j P_{(i,j)}(\tau^* \leq n-1)\right].$$

Notice that $\pi_{ij}^* = \pi_i \pi_j$ for $(i,j) \in \mathbb{Z}^2$ is an invariant measure for $(Y_n^*)$, and let $m^*$ be the shift-invariant measure defined on $(\mathbb{Z}^2)^\mathbb{Z}$ using $p^*$ and $\pi^*$ in the same way as measure $m$ was defined above using $p$ and $\pi$. Then

$$r_n = b_n^{-2\alpha}[P_{(0,0)}(\tau^* \geq n) + m^*(\tau^* \leq n-1)].$$

It follows from (5.2) that

$$P^*_{(0,0)}(Y^*_{nd} = (0,0)) = (P_0(x_{nd} = 0))^2 = n^{-2\gamma}(L(n))^2 \qquad \text{as } n \to \infty.$$

Therefore, for $\gamma \leq 1/2$ we may, once again, appeal to Remark 3.1 and Lemma 3.3 in Resnick, Samorodnitsky and Xue (2000) to conclude that

$$m^*(\tau^* \leq n-1) \sim \frac{d^{1-2\gamma}}{\Gamma(1+2\gamma)\Gamma(1-2\gamma)} n^{2\gamma}(L(n))^{-2}$$

as $n \to \infty$ and, hence,

$$r_n \to d^{-1} \frac{\Gamma(1+\gamma)^2 \Gamma(1-\gamma)^2}{\Gamma(1+2\gamma)\Gamma(1-2\gamma)}$$

as $n \to \infty$. In particular, the condition (4.8) holds in the case $\gamma = 1/2$, and so (5.3) holds in this case as well.

On the other hand, (4.8) fails in the case $0 < \gamma < 1/2$, and we conjecture that (5.3) does not hold either (cf. with Remark 4.4).

In the previous example we saw classes of stationary S$\alpha$S processes generated by conservative flows satisfying (4.6) and (4.8) and, hence, also (4.9). For these processes $b_n$ is regularly varying at infinity with exponent $\gamma \geq 1/(2\alpha)$. The following example exhibits stationary S$\alpha$S processes generated by conservative flows and also satisfying both (4.6) and (4.8), for which $b_n$ can be regularly varying at infinity with any exponent $0 < \gamma < 1/\alpha$.

EXAMPLE 5.4. Once again, let $\mathbf{X}$ be given by (5.1) with $d_\alpha = 1$, but now the control measure $m$ of the S$\alpha$S random measure $M$ is a probability measure under which the projections $g_k$, $k = \ldots, -1, 0, 1, 2, \ldots$ are i.i.d. positive Pareto random variables with

$$m(g_0 > x) = x^{-\theta} \qquad \text{for } x \geq 1$$

for some $\theta > \alpha$. Note that, for every $0 < p < \theta$,

(5.4) $$\int_E \max_{k=0,1,\ldots,n-1} g_k^p \, m(d\mathbf{g}) \sim c_{p,\theta} \, n^{p/\theta} \qquad \text{as } n \to \infty$$

for some finite positive constant $c_{p,\theta}$. Using (5.4) with $p = \alpha$ shows that

$$b_n \sim c_{\alpha,\theta}^{1/\alpha} n^{1/\theta} \qquad \text{as } n \to \infty,$$

and so (4.6) holds. Using, furthermore, (5.4) with $\alpha < p < \theta$ shows uniform integrability of the sequence $(b_n^{-\alpha} \max_{j=0,\ldots,n-1} |g_j(g)|^\alpha)$, $n = 1, 2, \ldots$, with respect to $m$. Since (4.11) is obvious in this case, we conclude that (4.8) holds as well.

That is, the convergence to extreme value distribution (4.9) holds here in all cases, and the exponent of regular variation $1/\theta$ of $b_n$ spans here the entire range $(0, 1/\alpha)$.

All of the examples thus far emphasize the message of Theorem 4.1: while the partial maxima grow at the rate of $n^{1/\alpha}$ for stationary S$\alpha$S processes with a nondegenerate part generated by a dissipative flow, the rate of growth is strictly slower if the process is generated by a conservative flow. In all these examples the actual rate of growth of $b_n$ and, hence, of the partial maximum was determined by the properties of the conservative flow. Our final example exhibits a new class of stationary S$\alpha$S processes generated by a conservative flow, in which the actual rate of growth of the partial maxima is determined by the kernel $f$ in the representation (2.2).

EXAMPLE 5.5. Let $E = (0,1]$, let $m$ be the Lebesgue measure on $E$, and define $\phi: E \to E$ by $\phi(x) = \{2x\}$, where $\{a\}$ is the fractional part of a real number $a$. Obviously, $(\phi^n, n = 0, 1, 2, \ldots)$ is a conservative measure-$m$-preserving flow (but it is not invertible). This flow is referred to as the *dyadic transformation*; see Example 1.6 in Billingsley (1965).

Let $h_1, h_2, \ldots$ be a nondecreasing sequence of nonnegative numbers such that

$$(5.5) \qquad \sum_{k=1}^{\infty} h_k^\alpha 2^{-k} < \infty,$$

and suppose that the limit

$$(5.6) \qquad \theta =: \lim_{k \to \infty} \frac{\log_2 h_k}{k} \in [0, 1/\alpha]$$

exists.

For $x = (x_1, x_2, \ldots) \in (0,1]$ given in its binary expansion, define

$$(5.7) \qquad K(x) = \inf\{j \geq 1 : x_j = 1\},$$

and let **X** be a stationary S$\alpha$S process defined by (2.1) with $f_n(x) = f \circ \phi^n(x)$, $x \in E$, where we choose

$$(5.8) \qquad f(x) = h_{K(x)}, \qquad x \in E.$$

Of course we do not need to worry about the fact that neither $K$ nor $f$ is defined for binary rational numbers. Notice, furthermore, that by switching to the space $E_1$ of doubly infinite sequences $x = (\ldots, x_{-1}, x_0, x_1, x_2, \ldots)$ with $x_i \in \{0,1\}$ for all $i$ with the product probability measure $m_1$ under which the projections $x_k$ form a sequence of i.i.d. Bernoulli random variables with mean $1/2$, and the left shift operator as the flow, we regain invertibility and, hence, we are in the framework of (2.2).

Denoting

$$(5.9) \qquad R_n(x) = \max_{j=1,\ldots,n} \sup\{m \geq 1 : x_j = x_{j+1} = \cdots = x_{j+m-1} = 0\}$$

to be the longest run of zeroes starting in $\{1,\ldots,n\}$, $n=1,2,\ldots$, we see that

(5.10) $$b_n^\alpha = \int_0^1 h_{R_n(x)}^\alpha \, dx, \qquad n=1,2,\ldots.$$

We claim that

(5.11) $$\lim_{k\to\infty} \frac{\log_2 b_k}{\log_2 k} = \theta.$$

Suppose that $\theta < 1/\alpha$. Then, for every $\epsilon \in (0, 1/\alpha - \theta)$, we have

$$b_n^\alpha \leq c_\epsilon \int_0^1 2^{\alpha(\theta+\epsilon)R_n(x)} \, dx, \qquad n=1,2,\ldots,$$

for some finite positive constant $c_\epsilon$. An elementary Bernoulli trials computation shows that, for any $0 < \gamma < 1$,

(5.12) $$\lim_{n\to\infty} \frac{\log_2 \int_0^1 2^{\gamma R_n(x)} \, dx}{\log_2 n} = \gamma.$$

Indeed, $R_n/\log_2 n \to 1$ a.s. [Example 4.14 in Billingsley (1986)], providing the liminf part of (5.12). Furthermore, for every $\delta > 0$ and $y > 0$,

$$\text{Leb}\{x : 2^{\gamma R_n(x)} > n^{\gamma(1+\delta)} + y\} \leq 2n^{-\delta} y^{-1/\gamma},$$

and so

$$\int_0^1 2^{\gamma R_n(x)} \, dx \leq (n^{\gamma(1+\delta)} + 1) + \int_1^\infty \text{Leb}\{x : 2^{\gamma R_n(x)} > n^{\gamma(1+\delta)} + y\} \, dy$$

$$\leq (n^{\gamma(1+\delta)} + 1) + 2n^{-\delta} \int_1^\infty y^{-1/\gamma} \, dy,$$

establishing the limsup part of (5.12).

We conclude that

$$\limsup_{k\to\infty} \frac{\log_2 b_k}{\log_2 k} \leq \theta + \epsilon,$$

and letting $\epsilon \to 0$ we see that, if $0 \leq \theta < 1/\alpha$, then

(5.13) $$\limsup_{k\to\infty} \frac{\log_2 b_k}{\log_2 k} \leq \theta.$$

Of course, the fact that (5.13) also holds for $\theta = 1/\alpha$ is trivial.

Similarly, if $\theta > 0$, then for every $\epsilon \in (0, \theta)$, we have

$$b_n^\alpha \geq c_\epsilon \int_0^1 2^{\alpha(\theta-\epsilon)R_n(x)} \, dx, \qquad n=1,2,\ldots,$$

for some finite positive constant $c_\epsilon$, and using, once again, (5.12) and letting $\epsilon \to 0$, we obtain

(5.14) $$\liminf_{k \to \infty} \frac{\log_2 b_k}{\log_2 k} \geq \theta$$

for $\theta > 0$. Since (5.14) is trivial for $\theta = 0$, we obtain (5.11) from (5.13) and (5.14).

One can easily see from Theorem 4.1 and a slight modification of the proof of (4.4) that (5.11) implies that

$$\frac{\log_2 M_n}{\log_2 n} \to \theta \quad \text{in probability}$$

as $n \to \infty$, which is a way of saying that $M_n$ grows as fast as $b_n$ [but it is not as precise as, say, (4.9)].

**6. Conclusions.** We have seen that the sequence of partial maxima of a stationary S$\alpha$S process generated by a dissipative flow grows always at the rate of $n^{1/\alpha}$ irrespective of the further properties of the flow or of the kernel in the integral representation of the process. On the other hand, the sequence of partial maxima of stationary S$\alpha$S processes generated by a conservative flow grows always at the rate slower than $n^{1/\alpha}$. However, the actual rate of growth here depends on the further properties of the flow and of the kernel. This is an example of a phase transition that would be consistent with a passage between short memory and long memory. It is important to note this phenomenon exists for all $0 < \alpha < 2$.

Examples of phase transitions of this kind for stationary S$\alpha$S processes have been observed before; see, for example, the discussion of ruin probabilities in Mikosch and Samorodnitsky (2000) or change in the rate of growth of partial sums in the case of increment processes of self-similar $\alpha$-stable processes with stationary increments; see Samorodnitsky and Taqqu (1994). In the previous cases the functionals for which the phase transitions were observed were based on partial sums of a process, which resulted in the phase transition being observed only in the range $1 < \alpha < 2$. In the present paper the phase transition occurs for all $0 < \alpha < 2$. It is interesting to mention that, unlike in the case of partial sums, maxima of long range dependent processes grow slower than those of short range dependent processes. This results in the fact that the maxima of a sum of two independent processes, one of which is generated by a dissipative flow, and the other by a conservative flow, will grow at the rate dictated by the dissipative part. If one decides to associate long range dependence of S$\alpha$S processes with conservative flows, then, for the partial maxima functional, the long memory is being hidden in such a sum. Recall, however, that this does not mean that the sum does

EXTREMES OF STATIONARY STABLE PROCESSES 31

not have long memory, for this cannot be dependent on behavior of a single functional.

It is also important to mention that, while stationary S$\alpha$S processes generated by conservative flows, considered in Mikosch and Samorodnitsky (2000), were also in the part of the parameter space suspected to represent long range dependence, some processes generated by dissipative flows were also long range dependent in that case.

This suggests that there may be multiple phase transitions indicating long range dependence, and the complete structure of those is yet to be discovered.

**Acknowledgment.** The author would like to thank the anonymous referee for the unusually careful reading of the paper and a number of useful detailed comments that helped the author to improve the paper.

School of Operations Research
and Industrial Engineering
Cornell University
Ithaca, New York 14853
USA
e-mail: gennady@orie.cornell.edu